\documentclass[10pt]{amsart}

\usepackage{amsfonts,amssymb}

\newtheorem{thm}{Theorem}[section]
\newtheorem{prop}[thm]{Proposition}

\newtheorem{lem}[thm]{Lemma}

\theoremstyle{remark}
\newtheorem{remark}{Remark}[section]

\theoremstyle{definition}

\newcommand*\isom{%
  \xrightarrow{\sim}%
}

\def\rr{\mathbb{R}}
\def\zz{\mathbb{Z}}

\def\cc{\mathbb{C}}
\def\pp{\mathbb{P}}

\def\hh{\mathbb{H}}

\def\uu{\Upsilon}

\def\isom{\, {\buildrel \sim \over \longrightarrow} \,}

\def\div{\mathrm{div} \, }

\def\deldelbar{\partial \overline{\partial}}

\begin{document}

\author{Robin de Jong}
\title{Gauss map on the theta divisor and Green's functions}



\begin{abstract}  In an earlier paper we constructed a 
Cartier divisor on the theta divisor of a 
principally polarised abelian variety whose support 
is precisely the ramification locus of the Gauss map. In this note we discuss a Green's function associated to this locus. For jacobians we relate
this Green's function to the canonical Green's function of the corresponding
Riemann surface.
\end{abstract}

\maketitle

\thispagestyle{empty}

\section{Introduction}

In \cite{jo} we investigated the properties of a certain theta function $\eta$
defined on the theta divisor of a principally polarised complex abelian variety
(ppav for short). Let us recall its definition. Fix a positive integer $g$ and
denote by $\hh_g$ the complex Siegel upper half space of degree $g$. On $\cc^g
\times \hh_g$ we have the Riemann theta function
\[ \theta = \theta(z,\tau) = \sum_{n \in \zz^g} \mathrm{e}^{\pi i {}^t n \tau n + 2\pi i
{}^t n z} \, . \]
Here and henceforth, vectors are column vectors and ${}^t$ denotes transpose.
For any fixed~$\tau$, the function $\theta = \theta(z)$ on $\cc^g$ gives rise to
an (ample, symmetric and reduced) divisor $\Theta$ on the torus
$A=\cc^g/(\zz^g + \tau \zz^g)$ which, by this token, acquires the structure of a
ppav. The theta function $\theta$ can be interpreted as a tautological section
of the line bundle $O_A(\Theta)$ on $A$.

Write $\theta_i$ for the first order partial derivative $\partial
\theta/\partial z_i$ and $\theta_{ij}$ for the second order partial derivative
$\partial^2 \theta / \partial z_i \partial z_j$.  
Then we define $\eta$ by 
\[ \eta = \eta(z,\tau) = \det \left( \begin{array}{cc} 
\theta_{ij} & \theta_j \\ {}^t \theta_i & 0 \end{array} \right)
\, . \]
We consider the restriction of $\eta$ to the vanishing locus of $\theta$ on
$\cc^g \times \hh_g$.

In \cite{jo} we proved that for any fixed $\tau$ the function $\eta$ gives rise
to a global section of the line bundle $O_\Theta(\Theta)^{\otimes g+1} \otimes
\lambda^{\otimes 2}$ on 
$\Theta$
in $A=\cc^g/(\zz^g+\tau \zz^g)$; here $\lambda$ is the trivial line bundle
$\mathrm{H}^0(A,\omega_A) \otimes_\cc O_\Theta$, with $\omega_A$ the canonical
line bundle on $A$. When viewed as a function of two
variables $(z,\tau)$ the function $\eta$ transforms as a theta function of
weight $(g+5)/2$ on $\theta^{-1}(0)$. If $\tau$ is fixed then the
support of $\eta$ on $\Theta$ is exactly the closure in
$\Theta$ of the ramification locus $\mathrm{R}(\gamma)$ of the Gauss map on the smooth
locus $\Theta^s$ of $\Theta$. Recall that the Gauss map on $\Theta^s$ 
is the map
\[ \gamma \colon \Theta^s \longrightarrow \pp(T_0 A)^\lor \]
sending a point $x$ in $\Theta^s$ to the tangent space $T_x \Theta$, translated
over $x$ to a subspace of $T_0 A$. It is well-known that the Gauss map on
$\Theta^s$ is generically finite exactly when $(A,\Theta)$ is
indecomposable; in particular the section $\eta$ is non-zero for such ppav's.

It turns out that the form $\eta$ has a rather nice application in the
study of the geometry of certain codimension-$2$ cycles on the moduli space of
ppav's. 
For this application we refer to the paper \cite{gsm}. 

The purpose of the present note is to discuss a certain real-valued variant
$\|\eta\|$ of $\eta$. 
In the case that $(A,\Theta)$ is the jacobian of a Riemann
surface $X$ we will establish a relation between this $\|\eta\|$ and the
canonical Green's function of $X$. In brief, note that in the case of a jacobian
of a Riemann surface $X$ we can
identify $\Theta^s$ with the set of effective 
divisors of degree $g-1$ on $X$ that do not
move in a linear system; thus for such divisors $D$ it makes sense to define
$\|\eta\|(D)$. On the other hand, note that $\Theta^s$ carries a canonical
involution $\sigma$ coming from the action of $-1$ on $A$, and moreover note
that sense can be made of evaluating the canonical (exponential) Green's
function $G$ of $X$ on pairs of effective divisors of $X$. The relation that we
shall prove is then of the form
\[ \|\eta\|(D) = \mathrm{e}^{-\zeta(D)} \cdot G(D,\sigma(D))  \, ; \]
here $D$ runs through the divisors in $\Theta^s$, and $\zeta$ is a certain
continuous function on $\Theta^s$. The $\zeta$ from the
above formula is
intimately connected with the geometry of intersections $\Theta \cap
(\Theta+R-S) $, where $R,S$ are distinct points on $X$. Amusingly, the limits of
such intersections where $R$ and $S$ approach each other 
are hyperplane sections of the Gauss
map corresponding to points on the canonical image of $X$, so the Gauss map on
the theta divisor is
connected with the above formula in at least two different ways. 
 
\section{Real-valued variant of $\eta$} \label{real}

Let $(A=\cc^g/(\zz^g + \tau \zz^g), \Theta = \div \theta)$ be a ppav as in the
introduction. As we said, the function $\eta$ transforms like a theta function
of weight $(g+5)/2$ and order $g+1$ on $\Theta$. This implies that if we define
\[ \|\eta\| = \|\eta\|(z,\tau) = (\det Y)^{(g+5)/4} \cdot 
\mathrm{e}^{-\pi (g+1) {}^t y \cdot Y^{-1} \cdot y} \cdot |\eta(z,\tau)| \, , \] 
where $Y= \mathrm{Im} \, \tau$ and $y = \mathrm{Im} \, z$, we obtain a
(real-valued) function which is invariant for the action of Igusa's transformation
group $\Gamma_{1,2}$ of matrices $\gamma = \left( \begin{array}{cc} a & b \\ c & d \end{array} \right)$ in
$\mathrm{Sp}(2g,\zz)$ with $a,b,c,d$ square matrices such that the diagonals of
both ${}^t a c$ and ${}^t b d$ consist of even integers. Recall that
$\Gamma_{1,2}$ acts on $\cc^g \times \hh_g$ via 
\[ (z,\tau) \mapsto \left( {}^t (c \tau + d)^{-1} z , (a\tau + b)(c \tau
+d)^{-1} \right) \, . \]
It follows that $\|\eta\|$ is a well-defined
function on $\Theta$, equivariant with respect to isomorphisms $(A,\Theta) \isom
(A',\Theta')$ coming from the symplectic action of $\Gamma_{1,2}$ on $\hh_g$.
Note that the zero locus of $\|\eta\|$ on $\Theta$ coincides with the zero locus of $\eta$
on $\Theta$. In fact, if $(A,\Theta)$ is indecomposable then 
the function $-\log \| \eta \|$ is a Green's function on
$\Theta$ associated to the closure of $\mathrm{R}(\gamma)$.

The definition of  $\|\eta\|$ is a variant upon the definition of the function
\[ \|\theta\| = \|\theta\|(z,\tau) = (\det Y)^{1/4} \cdot \mathrm{e}^{-\pi {}^t y \cdot
Y^{-1} \cdot y} \cdot |\theta(z,\tau)| \]
that one finds in \cite{fa}, p.~401. We note that $\|\theta\|$ should be seen as the
norm of $\theta$ for a canonical hermitian metric $\| \cdot \|_\mathrm{Th}$ on
$O_A(\Theta)$; we obtain $\|\eta\|$ as the norm of $\eta$ for the induced metric
on $O_\Theta(\Theta)^{\otimes g+1} \otimes \lambda^{\otimes 2}$. Here
$\mathrm{H}^0(A,\omega_A)$ has the standard metric given by putting $\|
\mathrm{d} z_1 \wedge \ldots \wedge \mathrm{d} z_g \| = (\det Y)^{1/2}$.  

The curvature form of $(O_A(\Theta), \| \cdot \|_\mathrm{Th})$ on $A$ is the
translation-invariant $(1,1)$-form
\[ \mu = \frac{\mathrm{i}}{2} \sum_{k=1}^g \mathrm{d} z_k \wedge \overline{\mathrm{d}z_k}  \, . \]
The $(g,g)$-form $\frac{1}{g!} \mu^g$ is a Haar measure for $A$ giving $A$
measure $1$. As $\mu$ represents $\Theta$ we have
\[ \frac{1}{g!} \int_\Theta \mu^{g-1} = 1 \, . \]
If $(A,\Theta)$ is indecomposable then $\log \|\eta\|$ is integrable with
respect to $\mu^{g-1}$ and the integral
\[ \frac{1}{g!} \int_\Theta \log \|\eta\| \cdot \mu^{g-1} \]
is a natural real-valued invariant of $(A,\Theta)$, which we think merits further attention.

\section{Arakelov theory of Riemann surfaces} \label{arakelov}

The purpose of this section and the next is 
to investigate the function $\|\eta\|$ in more detail for 
jacobians. There turns out to be a natural connection with certain real-valued
invariants occurring in the Arakelov theory of Riemann surfaces. We begin by
recalling the basic notions from this theory \cite{ar} \cite{fa}.

Let $X$ be a compact and connected Riemann surface of positive genus $g$, fixed
from now on. Denote
by $\omega_X$ its canonical line bundle. On $\mathrm{H}^0(X,\omega_X)$ we 
have a natural
inner product $(\omega,\eta) \mapsto \frac{\mathrm{i}}{2} \int_X \omega \wedge
\overline{\eta}$; we fix an orthonormal basis $(\omega_1,\ldots,\omega_g)$ with
respect to this inner product.

We put
\[ \nu = \frac{\mathrm{i}}{2g} \sum_{k=1}^g \omega_k \wedge \overline{\omega_k} \, . \]
This is a $(1,1)$-form on $X$, independent of our choice of 
$(\omega_1,\ldots,\omega_g)$ and hence 
canonical. In fact, if one denotes by $(J,\Theta)$ the jacobian of $X$
and by $j \colon X \hookrightarrow J$ an embedding of $X$ into $J$ using line integration, then
$\nu = \frac{1}{g} j^* \mu$ where $\mu$ is the 
translation-invariant form on $J$ discussed in the previous section. We have
$\int_X \nu =1$.
 
The canonical Green's function $G$ of $X$ is the unique non-negative 
function on $ X \times X$ which is non-zero outside the diagonal and satisfies
\[ \frac{1}{\mathrm{i} \pi} \partial \overline{\partial} \log G(P, \cdot) = \nu(P) -
\delta_P \, , \quad \int_X \log G(P,Q) \nu(Q) = 0  \]
for each $P$ on $X$; here $\delta$ denotes Dirac measure. The
functions $G(P,\cdot)$ give rise to canonical hermitian
metrics on the line bundles $O_X(P)$, with curvature form equal to $\nu$. 

From $G$, a smooth hermitian metric $\|\cdot \|_\mathrm{Ar}$ can be put on
$\omega_X$ by declaring that for each $P$ on $X$, the residue isomorphism
\[ \omega_X(P) [P] = (\omega_X \otimes_{O_X} O_X(P))[P] \isom \cc \]
is an isometry. Concretely this means that if $z : U \to \cc$ is a local
coordinate around $P$ on $X$ then 
\[ \|\mathrm{d} z\|_\mathrm{Ar}(P) = \lim_{ Q \to P} |z(P)-z(Q)|/G(P,Q) \, . \]
The curvature form of the metric $\|\cdot \|_\mathrm{Ar}$ on $\omega_X$ is equal
to $(2g-2) \nu$.

We conclude with the delta-invariant of $X$. 
Write $J=\cc^g/(\zz^g + \tau \zz^g)$ and $
\Theta = \div \theta$. There is a standard and
canonical identification of $(J,\Theta)$ with $(\mathrm{Pic}_{g-1} X, \Theta_0)$
where $ \mathrm{Pic}_{g-1} X$ is the set of linear equivalence classes of
divisors of degree $g-1$ on $X$, and where $\Theta_0 \subseteq
\mathrm{Pic}_{g-1} X$ is the subset of $\mathrm{Pic}_{g-1} X$ consisting of
the classes of effective divisors. By the identification 
$(J,\Theta) \cong ( \mathrm{Pic}_{g-1} X, \Theta_0)$ the function $\|\theta\|$
can be interpreted as a function on  $ \mathrm{Pic}_{g-1} X $. 

Now recall that the curvature form of $(O_J(\Theta), \| \cdot \|_\mathrm{Th})$ is
equal to $\mu$. This boils down to an equality of currents
\[ \frac{1}{\mathrm{i} \pi}  \deldelbar \log \|\theta\| = \mu - \delta_\Theta    \]
on $J$.
On the other hand one has for generic $P_1,\ldots,P_g$ on $X$
that $\|\theta\|(P_1+\cdots +P_g - Q)$ vanishes precisely when $Q$ is one of the
points $P_k$. This implies that on $X$ the equality of currents
\[ \frac{1}{\mathrm{i} \pi} \partial_Q \overline{\partial}_Q \log
\|\theta\|(P_1+\cdots+P_g-Q) = j^*\mu - \sum_{k=1}^g \delta_{P_k} = g \nu
- \sum_{k=1}^g \delta_{P_k}   \]
holds. Since also
\[ \frac{1}{\mathrm{i} \pi} \partial_Q \overline{\partial}_Q \log \prod_{k=1}^g G(P_k,Q)
= g  \nu - \sum_{k=1}^g \delta_{P_k} \]
we may conclude, by compactness of $X$, that
\[ \|\theta\|(P_1+\cdots+P_g-Q) = c(P_1,\ldots,P_g) \cdot \prod_{k=1}^g G(P_k,Q)
\] for some constant $c(P_1,\ldots,P_g)$ depending only on $P_1,\ldots,P_g$. A
closer analysis (cf. \cite{fa}, p.~402) reveals that
\[ c(P_1,\ldots,P_g) = \mathrm{e}^{-\delta/8} \cdot \frac{ \| \det \omega_i (P_j)
\|_\mathrm{Ar} }{ \prod_{k<l} G(P_k,P_l) } \]
for some constant $\delta$ which is then by definition the delta-invariant of
$X$. The argument to prove this equality 
uses certain metrised line bundles and their
curvature forms on sufficiently big powers $X^r$ of $X$. A variant of this
argument occurs in the proof of our main result below.
  
\section{Main result}

In order to state our result, we need some more notation and facts. 
We still have our fixed Riemann surface $X$ of positive genus $g$ 
and its jacobian $(J,\Theta)$. 
The following lemma is well-known.
\begin{lem} Under the identification $\Theta \cong \Theta_0$, the smooth locus
$\Theta^s$ of $\Theta$ corresponds to the subset $\Theta_0^s$ of $\Theta_0$ of
divisors that do not move in a linear system. Furthermore, there is a
tautological surjection $\Sigma$ from the $(g-1)$-fold symmetric power
$X^{(g-1)}$ of $X$ onto $\Theta_0$. This map $\Sigma$ is an isomorphism over
$\Theta_0^s$.
\end{lem}
The lemma gives rise to identifications $\Theta^s \cong
\Theta_0^s \cong \uu$ with $\uu$ a certain open subset of $X^{(g-1)}$. 
We fix and accept these identifications in all that follows.
Note that the set $\uu$ carries a
canonical involution $\sigma$, coming from the action of $-1$ on $J$. 
For $D$ in $\uu$ the divisor $D + \sigma(D)$ of degree $2g-2$ is always 
a canonical divisor. 

The next lemma gives a description of the ramification locus of the Gauss map 
on $\Theta^s \cong \uu$.
\begin{lem} Under the identification $\Theta^s \cong \uu$ the ramification
locus of the Gauss map on $\Theta^s$ corresponds to the set of divisors $D$ in $\uu$ such that
$D$ and $\sigma(D)$ have a point in common.
\end{lem}   
\begin{proof} According to \cite{deb}, p.~691 the ramification locus of the
Gauss map is given by the set of divisors $E+P$ with $E$ effective of degree
$g-2$ and $P$ a point such that on the canonical image of $X$ the divisor 
$E+2P$ is contained in a hyperplane. 
But this condition on $E$ and $P$ means that $E+2P$ is
dominated by a canonical divisor, or equivalently, that 
$P$ is contained in the conjugate $\sigma(E+P)$ of $E+P$. The lemma follows. 
\end{proof}

If $D=P_1+\cdots+P_m$ and $D'=Q_1+\cdots+Q_n$ are two effective divisors on $X$
we define $G(D,D')$ to be
\[ G(D,D') = \prod_{i=1}^m \prod_{j=1}^n G(P_i,Q_j) \, . \]
Clearly the value $G(D,D')$ is zero if and only
if $D$ and $D'$ have a point in common. Applying this to the above lemma, we see that
the function $D \mapsto G(D,\sigma(D))$ on $\uu$ vanishes precisely on the
ramification locus of the Gauss map. As a consequence $G(D,\sigma(D))$ and
$\|\eta\|(D)$ have exactly the same zero locus. It looks therefore as if 
a relation 
\[ \|\eta\|(D) = \mathrm{e}^{-\zeta(D)} \cdot G(D,\sigma(D)) \]
should hold for $D$ on $\uu$ with $\zeta$ a suitable continuous function. 
The aim of the rest of this note is to prove 
this relation, and to compute $\zeta$ explicitly. 

We start with
\begin{prop} Let $Y=\uu \times X \times X $.
The map $\| \Lambda \| \colon Y \to \rr$ given
by
\[ \|\Lambda\|(D,R,S) = \frac{ \| \theta \| (D+R-S) }{ G(R,S) G(D,S)
G(\sigma(D),R) } \]
is continuous and nowhere vanishing. 
Furthermore $\|\Lambda \|$ factors via the 
projection of $Y$ onto $\uu$.
\end{prop}
\begin{proof} The numerator $\|\theta\|(D+R-S)$ vanishes if and only if $R=S$ or
$D=E+S$ for some effective divisor $E$ of degree $g-2$ or $D+R$ is linearly
equivalent to an effective divisor $E'$ of degree $g$ such that $E'=E''+S$ for
some effective divisor $E''$ of degree $g-1$. The latter condition is precisely
fulfilled when the linear system $|D+R|$ is positive dimensional, or
equivalently, by Riemann-Roch, when $D+R$ is dominated by a canonical divisor,
i.e. when $R$ is contained in $\sigma(D)$. It follows that the
numerator  $\|\theta\|(D+R-S)$ and the denominator $ G(R,S) G(D,S)
G(\sigma(D),R)$ have the same zero locus on $Y$. Fixing a divisor 
$D$ in $\uu$ and using
what we have said in Section \ref{arakelov} it is seen that the currents
\[ \frac{1}{\mathrm{i} \pi} \deldelbar \log \|\theta\|(D+R-S) \,\, \textrm{and}
\,\,
\frac{1}{\mathrm{i} \pi} \deldelbar \log \left( G(R,S) G(D,S)
G(\sigma(D),R) \right) \]
are both the same on $X \times X$. We conclude that $\|\Lambda \|$ is non-zero
and continuous and depends only on $D$.
\end{proof}
We also write $\|\Lambda\|$ for the induced map on $\uu$.
Our main result is
\begin{thm} \label{main}
Let $D$ be an effective divisor of degree $g-1$ on $X$, not moving in a linear
system. Then the formula
\[ \| \eta \|(D) = \mathrm{e}^{-\delta/4} \cdot \|\Lambda\|(D)^{g-1} \cdot G(D,\sigma(D)) \]
holds.
\end{thm} 
\begin{proof} Fix two distinct points $R,S$ on $X$. We start
by proving that there is a non-zero constant $c$ depending only on $X$ such that
\[ (*) \quad \|\eta\|(D) = c \cdot G(D,\sigma(D)) 
\left( \frac{ \| \theta \| (D+R-S) }{ G(R,S) G(D,S) G(\sigma(D),R)}
\right)^{g-1}  \]  
for all $D$ varying through $\uu$. We would be done if we could show that
\[ \frac{1}{\mathrm{i} \pi} \deldelbar \log \|\eta\|(D) \]
and
\[ \frac{1}{\mathrm{i} \pi} \deldelbar \log \left( G(D,\sigma(D)) 
\left( \frac{ \| \theta \| (D+R-S) }{ G(R,S) G(D,S) G(\sigma(D),R)}
\right)^{g-1} \right) \]
define the same currents on $\uu$. Indeed, then the function $\phi(D)$ given by
\[ \log \|\eta\|(D) - \log \left( G(D,\sigma(D)) 
\left( \frac{ \| \theta \| (D+R-S) }{ G(R,S) G(D,S) G(\sigma(D),R)}
\right)^{g-1} \right) \]
is pluriharmonic on $\uu$, hence on $\Theta^s$, and since $\Theta^s$ is open in
$\Theta$ with boundary empty or of codimension $\geq 2$, and since $\Theta$ is
normal (cf. \cite{ke}, Theorem~1') we may conclude that $\phi$ is constant.

To prove equality of
\[ \frac{1}{\mathrm{i} \pi} \deldelbar \log \|\eta\|(D) \]
and
\[ 
\frac{1}{\mathrm{i} \pi} \deldelbar \log \left( G(D,\sigma(D)) 
\left( \frac{ \| \theta \| (D+R-S) }{ G(R,S) G(D,S) G(\sigma(D),R)}
\right)^{g-1} \right) \]
on $\uu$ it suffices to prove that their pullbacks are equal on $\uu' =
p^{-1}(\uu) $ in $ X^{g-1}$ under the canonical projection $p : X^{g-1} \to
X^{(g-1)}$. 

First of all we compute the pullback under $p$ of
\[ \frac{1}{\mathrm{i} \pi} \deldelbar \log \|\eta\|(D) \]
on $\uu'$. Let $\pi_i \colon X^{g-1} \to X$ for $i=1,\ldots,g-1$ be the
projections onto the various factors. We have seen that the curvature form of
$O_J(\Theta)$ is $\mu$, hence the curvature form of $O_\Theta(\Theta)^{\otimes
g+1}$ is $(g+1)\mu_\Theta$. According to \cite{fa}, p.~397 the pullback of
$\mu_\Theta$ to $X^{g-1}$ under the canonical surjection 
$\Sigma \colon X^{g-1} \to \Theta$ can be written as
\[ \frac{\mathrm{i}}{2} \sum_{k=1}^g \left(  \sum_{i=1}^{g-1} \pi_i^*(\omega_k) \right) \wedge 
\left(\sum_{i=1}^{g-1} \pi_i^* (\overline{\omega_k}) \right) \, . \]
Here $(\omega_1,\ldots,\omega_g)$ is an orthonormal basis for $\mathrm{H}^0(X,\omega_X)$
which we fix. Let's call the above form $\xi$. It follows that
\[ p^* \frac{1}{\mathrm{i} \pi } \deldelbar \log \|\eta\|(D) = (g+1) \xi -
\delta_{p^* \mathrm{R}(\gamma) } \]
as currents on $\uu'$. Here $\mathrm{R}(\gamma)$ is the ramification locus of the Gauss
map on $\uu$.

Next we consider the pullback under $p$ of
\[ \frac{1}{\mathrm{i} \pi} \deldelbar \log \left( G(D,\sigma(D)) 
\left( \frac{ \| \theta \| (D+R-S) }{ G(R,S) G(D,S) G(\sigma(D),R)}
\right)^{g-1}  \right)   \]
on $\uu'$. The factor $\|\theta\|(D+R-S)$ accounts for a contribution equal
to $\xi$, and both of the factors $G(D,S)$ and $G(\sigma(D),R)$ give a
contribution $\sum_{i=1}^{g-1} \pi_i^*(\nu)$. We find
\[ p^* \frac{1}{\mathrm{i} \pi } \deldelbar \log \left( \frac{ \| \theta \| (D+R-S) }{ G(R,S) G(D,S) G(\sigma(D),R)}
\right)^{g-1} = (g-1)  ( \xi - 2 \sum_{i=1}^{g-1} \pi_i^* (\nu)  ) \, . \]
We are done if we can prove that
\[ p^* \frac{1}{\mathrm{i} \pi } \deldelbar \log G(D,\sigma(D)) = 2\xi + (2g-2)
\sum_{i=1}^{g-1} \pi_i^*(\nu) - \delta_{p^* \mathrm{R}(\gamma)} \, . \]
For this consider the product $\uu' \times \uu' \subseteq X^{g-1} \times
X^{g-1}$. For $i,j=1,\ldots,g-1$ denote by $\pi_{ij} : X^{g-1} \times X^{g-1}
\to X \times X$ the projection onto the $i$-th factor of the left $X^{g-1}$, and
onto the $j$-th factor of the right $X^{g-1}$. Denoting by $\Phi$ the smooth
form represented by $\frac{1}{\mathrm{i} \pi } \deldelbar \log G(P,Q)$ on $X \times X$ it
is easily seen that
we can write
\[ p^* \frac{1}{\mathrm{i} \pi } \deldelbar \log G(D,\sigma(D)) + 
\delta_{p^* \mathrm{R}(\gamma)}
= ( \sigma^* \sum_{i,j=1}^{g-1} \pi_{ij}^* \Phi )\arrowvert_\Delta \, ; \]
here $\Delta \cong \uu'$ is the diagonal in $\uu' \times \uu'$ and $\sigma^*$ is
the action on symmetric $(1,1)$-forms on $\uu'$ induced 
by the automorphism $(x,y) \mapsto (x,\sigma(y))$ of $\uu \times \uu$. 
Let $q_1,q_2$ be the projections of $X \times X$ onto the first and second
factor, respectively. Then according
to \cite{ar}, Proposition 3.1 we have
\begin{eqnarray*} \Phi & = & \frac{\mathrm{i}}{2g} \sum_{k=1}^g q_1^*(\omega_k) \wedge
q_1^*(\overline{\omega_k}) +  \frac{\mathrm{i}}{2g} \sum_{k=1}^g q_2^*(\omega_k) \wedge
q_2^*(\overline{\omega_k}) \\ & & - \frac{\mathrm{i}}{2} \sum_{k=1}^g q_1^*(\omega_k) \wedge
q_2^*(\overline{\omega_k}) - \frac{\mathrm{i}}{2} \sum_{k=1}^g q_2^*(\omega_k) \wedge
q_1^*(\overline{\omega_k}) \, .
\end{eqnarray*}
Note that $q_1 \cdot \pi_{ij} = \pi_i$ and $q_2 \cdot
\pi_{ij} = \pi_j$; this gives
\begin{eqnarray*} \pi_{ij}^* \Phi & = & \frac{\mathrm{i}}{2g} \sum_{k=1}^g \pi_i^*(\omega_k) \wedge
\pi_i^*(\overline{\omega_k}) +  \frac{\mathrm{i}}{2g} \sum_{k=1}^g \pi_j^*(\omega_k) 
\wedge
\pi_j^*(\overline{\omega_k}) \\ & & - \frac{\mathrm{i}}{2} \sum_{k=1}^g \pi_i^*(\omega_k) \wedge
\pi_j^*(\overline{\omega_k}) - \frac{\mathrm{i}}{2} \sum_{k=1}^g \pi_j^*(\omega_k) \wedge
\pi_i^*(\overline{\omega_k}) \, .
\end{eqnarray*}
Next note that $\sigma$ acts as $-1$ on $\mathrm{H}^0(X,\omega_X)$; this implies, at least
formally, that
\begin{eqnarray*} \sigma^* \pi_{ij}^* \Phi  & = & \frac{\mathrm{i}}{2g} \sum_{k=1}^g \pi_i^*(\omega_k) \wedge
\pi_i^*(\overline{\omega_k}) +  \frac{\mathrm{i}}{2g} \sum_{k=1}^g \pi_j^*(\omega_k) 
\wedge
\pi_j^*(\overline{\omega_k}) \\ & & + \frac{\mathrm{i}}{2} \sum_{k=1}^g \pi_i^*(\omega_k) \wedge
\pi_j^*(\overline{\omega_k}) + \frac{\mathrm{i}}{2} \sum_{k=1}^g \pi_j^*(\omega_k) \wedge
\pi_i^*(\overline{\omega_k}) \\ & = & \pi_i^* (\nu) + \pi_j^* (\nu) + 
\frac{\mathrm{i}}{2} \sum_{k=1}^g \pi_i^*(\omega_k) \wedge
\pi_j^*(\overline{\omega_k}) \\ & & + \frac{\mathrm{i}}{2} \sum_{k=1}^g \pi_j^*(\omega_k) \wedge
\pi_i^*(\overline{\omega_k}) \, .
\end{eqnarray*}
We obtain for $ ( \sigma^* \sum_{i,j=1}^{g-1} \pi_{ij}^* \Phi
)\arrowvert_\Delta$ the expression
\begin{eqnarray*}     
\sum_{i,j=1}^{g-1}  \pi_i^* (\nu) + \pi_j^* (\nu)     +
\frac{\mathrm{i}}{2} \sum_{k=1}^g \sum_{i,j=1}^{g-1}    
\pi_i^*(\omega_k) \wedge
\pi_j^*(\overline{\omega_k}) + \pi_j^*(\omega_k) \wedge
\pi_i^*(\overline{\omega_k})   \\
  =   (2g-2) \sum_{i=1}^{g-1} \pi_i^*(\nu) + \mathrm{i} \sum_{k=1}^g  \left(
(\sum_{i=1}^{g-1} \pi_i^* (\omega_k)) \wedge ( \sum_{i=1}^{g-1}
\pi_i^*(\overline{\omega_k}))  \right) \\   =  (2g-2) \sum_{i=1}^{g-1}
\pi_i^*(\nu) + 2 \xi \, , \hspace{4.97cm}
\end{eqnarray*}
and this gives us what we want. 

It remains to prove that the constant $c$ is equal to $\mathrm{e}^{-\delta/4}$.
We use the following lemma. 
\begin{lem} \label{second} Let $\mathrm{Wr}(\omega_1,\ldots,\omega_g)$ be the Wronskian
on $(\omega_1,\ldots,\omega_g)$, considered as a global section of
$\omega_X^{\otimes g(g+1)/2}$. Let $P$ be any point on $X$.
Then the equality
\[ \|\eta\| ((g-1)  P) = \mathrm{e}^{-(g+1)\delta/8} \cdot \|
\mathrm{Wr}(\omega_1,\ldots,\omega_g) \|_\mathrm{Ar}(P)^{g-1} \]
holds. Left and right hand side are non-vanishing for generic $P$.
\end{lem}
\begin{proof} 
Let $\kappa \colon X \to \Theta$ be the map given by sending $P$ on $X$ to the
linear equivalence class of $(g-1) \cdot P$.
According to \cite{jo2}, Lemma~3.2 we have a canonical isomorphism
\[ \kappa^*(O_\Theta(\Theta)) \otimes \omega_X^{\otimes g} \isom \omega_X^{\otimes
g(g+1)/2} \otimes \kappa^*(\lambda)^{\otimes -1} \]
of norm $\mathrm{e}^{\delta/8}$. It follows that we have a canonical isomorphism
\[
\kappa^*\left(O_\Theta(\Theta)^{\otimes g+1} \otimes \lambda^{\otimes 2} \right) 
\isom  \left( \omega_X^{\otimes g(g+1)/2} \otimes \kappa^*(\lambda)^{\otimes -1}
\right)^{\otimes g-1}
\]
of norm $\mathrm{e}^{(g+1)\delta/8}$. Chasing these isomorphisms using
\cite{jo}, Theorem~5.1 one sees that the global section $\kappa^* \eta$ of 
\[ \kappa^*\left(O_\Theta(\Theta)^{\otimes g+1} \otimes \lambda^{\otimes 2} \right) \] 
is sent to the global section
\[ \left( \xi_1 \wedge \ldots \wedge \xi_g \mapsto \frac{\xi_1 \wedge \ldots
\wedge \xi_g}{\omega_1 \wedge \ldots \wedge \omega_g} \cdot
\mathrm{Wr}(\omega_1,\ldots,\omega_g)  \right)^{\otimes g-1}  \]
of 
\[ \left( \omega_X^{\otimes g(g+1)/2} \otimes \kappa^*(\lambda)^{\otimes -1}
\right)^{\otimes g-1} \, . \]
The claimed equality follows.
The non-vanishing for generic
$P$ follows from the Wronskian being
non-zero as a section of $\omega_X^{\otimes g(g+1)/2}$. 
\end{proof}
We can now finish the proof of Theorem \ref{main}. Using the defining relation
\[ \|\theta\|(P_1+\cdots+P_g-S) = \mathrm{e}^{-\delta/8} 
\frac{ \| \det \omega_i (P_j)
\|_\mathrm{Ar} }{ \prod_{k<l} G(P_k,P_l) } \prod_{k=1}^g G(P_k,S) \]
mentioned earlier for $\delta$ we can rewrite equality (*) as
\[ \|\eta\|(D)  = c \cdot  \mathrm{e}^{-(g-1)\delta/8}   \frac{
G(D,\sigma(D))}{G(R,\sigma(D))^{g-1}}  
\left( \frac{ \| \det \omega_i (P_j)
\|_\mathrm{Ar} }{ \prod_{k<l} G(P_k,P_l) } \right)^{g-1} \, ; \]
here we have set $D=P_1+\cdots+P_{g-1}$ and 
$P_g=R$. Letting the $P_j$ approach $R$ we find, by a similar computation as in
\cite{jo2}, proof of Lemma~3.2,
\[ \|\eta\|((g-1)R) = c \cdot  \mathrm{e}^{-(g-1)\delta/8} \cdot
\|\mathrm{Wr}(\omega_1,\ldots,\omega_g) \|_\mathrm{Ar}(R)^{g-1} \, . \]
Lemma \ref{second} gives $c \cdot
\mathrm{e}^{-(g-1)\delta/8}=\mathrm{e}^{-(g+1)\delta/8}$, in other words
$c=\mathrm{e}^{-\delta/4}$.
\end{proof}

\begin{remark} \label{rem} 
It was shown by J.-B. Bost \cite{bo} that there is an invariant $A$
of $X$ such that for each pair of distinct points $R,S$ on $X$ the formula
\[ \log G(R,S) = \frac{1}{g!} \int_{\Theta+R-S} \log \|\theta\| \cdot \mu^{g-1}
+ A \]
holds. An inspection of the proof as for example given in \cite{we}, Section~5
reveals that the integrals
\[ \frac{1}{g!} \int_{\Theta^s} \log G(D,S) \cdot \mu(D)^{g-1} \quad
\textrm{and} \quad \frac{1}{g!} \int_{\Theta^s} \log G(\sigma(D), R)
\cdot \mu(D)^{g-1} \]
are zero and hence from the definition of $\|\Lambda\|$ we can write
\[ A = -\frac{1}{g!} \int_{\Theta^s} \log \|\Lambda\|(D) \cdot \mu(D)^{g-1}   \, .
\]

\end{remark}

\subsection*{Acknowledgements} The author is supported by a grant from the Netherlands Organisation for Scientific Research (NWO).

\vspace{0.5cm}

\noindent Address of the author:\\ \\
Mathematical Institute \\
University of Leiden \\
PO Box 9512 \\
2300 RA Leiden \\
The Netherlands \\
Email: \verb+rdejong@math.leidenuniv.nl+

\end{document}